\documentclass[10pt,a4paper]{amsart}
\usepackage{amsmath, amsthm, amsfonts, amssymb, amstext}

\newtheorem{theorem}{Theorem}[section]
\newtheorem{proposition}[theorem]{Proposition}
\newtheorem{lemma}[theorem]{Lemma}
\newtheorem{corollary}[theorem]{Corollary}

\theoremstyle{definition}

\newtheorem{remark}[theorem]{Remark}
\newtheorem{example}[theorem]{Example}

\begin{document}
\title[Automorphism Groups of Semifields]
{Automorphism Groups of some Finite Semifields}
\author{Andrew Steele}
\address{School of Mathematical Sciences, University of Nottingham, University Park, Nottingham, NG7 2RD, United Kingdom}
\email{pmxas4@nottingham.ac.uk}


\begin{abstract}
We determine the automorphism group for some well known constructions of finite semifields. In particular, we compute the automorphism group of Sandler's semifields and in certain cases the automorphism groups of the Hughes-Kleinfeld and Knuth semifields. We also determine how many nonisomorphic Sandler semifields can be constructed given a finite field $F$ and a finite extension $L/F$.
\end{abstract}
\maketitle
\section*{Introduction}
Finite semifields are finite, nonassociative, unital division algebras. They are traditionally studied in the context of finite geometries due to their connections with projective planes. In fact, every proper semifield coordinatizes a non-Desarguesian projective plane and two semifields coordinatize the same projective plane if and only if they are \emph{isotopic} \cite{MR0116036}. Because of this, semifields are usually classified up to isotopy rather than up to isomorphism and in many cases the automorphism groups of the semifields are not known. Finite semifields have also found applications in coding theory \cite{MR1455862}, \cite{MR1984461}, \cite{MR2287192} and combinatorics and graph theory \cite{MR2402272}. In this paper we study them as algebraic objects in their own right. The methods used are purely algebraic, in particular, we classify Sandler's semifields up to isomorphism and we study the automorphism group of this and several other classes of semifields. The main constructions we consider are the Sandler semifields \cite{MR0147519}, the Hughes-Kleinfeld semifields \cite{MR0117261} and the Knuth semifields \cite{MR0175942}.

The structure of the paper is as follows: in Section 1 basic preliminaries and notations used are introduced. In Section 2 we outline the construction of Sandler's semifields; in fact, we give a more general construction over infinite base fields following \cite{NACA} and recall any results necessary for what follows. In Section 3 we study Sandler's semifields up to isomorphism and compute their automorphisms. In Sections 4 and 5 we recall the definition of the Hughes-Kleinfeld and the Knuth semifields and calculate their nuclei which is the key result used in Section 7 to calculate the automorphisms of these semifields. This solves one of the open problems mentioned in Section 5 of \cite{MR2942770}.

This work will form part of the authors PhD thesis written under the supervision of Dr S. Pumpl\"{u}n.

\section{Preliminaries}
A \emph{finite semifield} is a finite set $S$ with at least two distinct elements. In addition, $S$ possesses two binary operations, $+$ and $\circ$, which satisfy the following axioms:
\begin{enumerate}
\item $(S,+)$ is a group with identity 0.
\item If $a$ and $b$ are elements of $S$ and $a \circ b=0$ then $a=0$ or $b=0$.
\item Distributivity holds: $a \circ (b+c) = a \circ b+a \circ c$ for all $a,b,c \in S$.
\item There is a multiplicative identity 1: $a \circ 1 = 1 \circ a = a$ for all $a \in S$.
\end{enumerate}
For ease of notation we will denote the multiplication operation $\circ$ by juxtaposition: $x \circ y = xy$ and throughout by `semifield' we mean `finite semifield'. It can be shown that semifields posses a vector space structure over some prime field $F=\mathbb{F}_p$ (see \cite{MR2939390} for example), so the number of elements in a finite semifield $S$ is $p^n$ where $n$ is the dimension of $S$ over $F$. Hence semifields are in fact (not necessarily associative) division algebras over finite fields. A famous theorem of Wedderburn \cite{MR1500717} tells us that every finite, associative division algebra is a finite field. Hence, any finite semifield is necessarily either a finite field or it is \emph{not} associative. A semifield which is not a finite field is called a \emph{proper semifield}. The number of elements in a semifield $S$ is called the \emph{order} of $S$. The \emph{associator} of three elements $x,y,z \in S$ is defined by
\[[x,y,z] = (xy)z - x(yz),\]
and is used to measure associativity in semifields. 

Semifields possess the following important substructures which are all isomorphic to a finite field:
\begin{align*}
Nuc_l(S) &:= \{x \in S \mid [x,y,z]=0 \text{ for all } y,z \in S\}, \\
Nuc_m(S) &:= \{y \in S \mid [x,y,z]=0 \text{ for all } x,z \in S\},\\
Nuc_r(S) &:= \{z \in S \mid [x,y,z]=0 \text{ for all } x,z \in S\}.
\end{align*}
These are called the \emph{left, middle} and \emph{right nuclei} respectively. The intersection of these three sets is called the \emph{nucleus} of $S$ denoted Nuc($S$). The \emph{commutative centre} of $S$ is the set
\[K(S) := \{x \in S \mid xy = yx \text{ for all } y \in S\}.\]
The intersection of the commutative centre and the nucleus of $S$ is called the \emph{centre} and is denoted $Z(S)$. For a survey of results on finite semifields we refer the reader to \cite{MR1725526} or \cite{MR2258004}. 

\section{Nonassociative Cyclic Algebras}
In this Section we outline Sandler's construction of a semifield of order $q^{n^2}$ were $q$ is a prime power and $n$ is an integer strictly bigger than $1$. This construction was studied by the author over arbitrary base fields in \cite{NACA} hence in this Section we follow the notation used in that paper and proofs of theorems can be found therein. In the above mentioned paper, these semifields are called \emph{nonassociative cyclic algebras}. This is due to the fact that they were originally constructed as a generalisation of Waterhouse's nonassociative quaternion algebras \cite{MR923618}. It was only later pointed out to the author that they were also generalisations of Sandler's semifields to infinite base fields.

Let $F$ be a field and let $L$ be a cyclic extension of $F$ of degree $n$, i.e. a Galois field extension with Galois group Gal$(L/F) = \mathbb{Z}/n\mathbb{Z}$. Pick an element $a \in L \setminus F$ and define the nonassociative cyclic algebra $(L/F, \sigma, a)$ to be the $L$-vector space 
\[(L/F, \sigma, a) := \bigoplus_{i=0}^{n-1} Lz^i\]
with basis $\{z^0 = 1, z, z^2, \ldots, z^{n-1}\}$ where the $z^i$ are formal symbols. We define a multiplication on elements $lz^i$ and $mz^j$ for $l,m \in L, 0 \leq i,j, <n$, by

\[
 (lz^i)(mz^j) =
  \begin{cases}
   l \sigma^i(m) z^{i+j} & \text{if } i+j < n \\
   l \sigma^i(m) a z^{(i+j)-n} & \text{if } i+j \geq n
  \end{cases}
\]
and then extend it linearly to all of $(L/F,\sigma,a)$. We say that the above algebra is a \emph{nonassociative cyclic algebra of degree $n$}.
It turns out that for many choices of $a \in L\setminus F$ this construction gives us a division algebra.

\begin{theorem}\label{div}
Let $(L/F, \sigma, a)$ be a nonassociative cyclic algebra of degree $n$. If the elements $1,a,a^2, \ldots, a^{n-1}$ are linearly independent over $F$ then $(L/F, \sigma, a)$ is a division algebra.
\end{theorem}

In particular, this theorem implies that if $L/F$ is a cyclic extension of prime degree then every choice of $a \in L \setminus F$ gives a division algebra of the form $(L/F, \sigma, a)$. 

\begin{theorem}\label{nuc}
Let $A$ be the nonassociative cyclic algebra $(L/F,\sigma, a)$ of degree $n$. Then
\begin{enumerate}
\item $Nuc(A) = L$.
\item $Z(A) = F$.
\end{enumerate}
\end{theorem}

We note here that the middle and right nuclei are also equal to the field extension $L$. However, the left nucleus depends on the choice of the element $a \in L \setminus F$. We make this explicit here since there appears to be a small mistake concerning this in Sandler's original paper.

Recall that by the definition of $(L/F, \sigma, a)$ we cannot have $a \in F$. However, $a$ may belong to some proper subfield $E \subset F$. In this case there will exist a proper subgroup $G_E$ of $Gal(L/F)$ such that for all $\tau \in G_E$ we have $\tau(a) = a$. Since $Gal(L/F)$ is a cyclic group, the subgroup $G_E$ will be generated by some power of the generator of $Gal(L/F)$. If the generator of $Gal(L/F)$ is $\sigma$ then denote the generator of $G_E$ by $\sigma^s$ where $2 \leq s \leq n-1$.

\begin{proposition}
Let $A = (L/F, \sigma, a)$ where $a$ belongs to a proper subalgebra $E$ of $L$. Then, using the notation of the previous paragraph,
\[Nuc_l(A) = L \oplus Lz^s \oplus Lz^{2s} \oplus \cdots \oplus Lz^{n-s},\]
where $s$ is such that $a$ is invariant under the subgroup $\langle \sigma^s \rangle$ of $Gal(L/F)$, i.e. $\sigma^{ks}(a) = a$ for all $k \in \mathbb{Z}$. 
\begin{proof}
Let 
\[B= L \oplus Lz^s \oplus Lz^{2s} \oplus \cdots \oplus Lz^{n-s},\]
which is clearly a subalgebra of $A$. By the distributivity of the multiplication in $A$ it suffices to check associativity for monomials in $A$. Let $wz^{m} \in B$ where $m$ is some multiple of $s$. Therefore $\sigma^m(a) = a$. Also let $xz^i, yz^j$ be monomials in  $A$ where $0 \leq i,j \leq n-1$. We have
\begin{align*}
wz^m((xz^i)(yz^j)) &= \begin{cases}
   wz^m(x \sigma^i(y) z^{i+j}) & \text{if } i+j < n \\
   wz^m(x \sigma^i(y) a z^{(i+j)-n}) & \text{if } i+j \geq n
  \end{cases}\\
&= \begin{cases}
   w \sigma^m(x) \sigma^{i+m}(y) z^{i+j+m} & \text{if } i+j+m < n \\
   w\sigma^m(x) \sigma^{i+m}(y) a z^{(i+j+m)-n} & \text{if } n \leq i+j+m < 2n\\
   w \sigma^m(x)\sigma^{i+m}(y)a^2 z^{(i+j+m)-2n} & \text{if } i + j + m \geq 2n
  \end{cases}
\end{align*}
whereas 
\begin{align*}
((wz^m)(xz^i))(yz^j) &= \begin{cases}
   (w\sigma^m(x)z^{m+i}) y z^{j} & \text{if } m+i < n \\
   (w\sigma^m(x)az^{(m+i)-n}) y  z^{j} & \text{if } m+i \geq n
  \end{cases}\\
&= \begin{cases}
   w \sigma^m(x) \sigma^{m+i}(y) z^{i+j+m} & \text{if } i+j+m < n \\
   w\sigma^m(x) \sigma^{m+i}(y) a z^{(i+j+m)-n} & \text{if } n \leq i+j+m < 2n\\
   w \sigma^m(x)\sigma^{m+i}(y)a^2 z^{(i+j+m)-2n} & \text{if } i + j + m \geq 2n.
  \end{cases}
\end{align*}
In any of the cases these terms are equal so the subalgebra mentioned in the proposition is contained in the left nucleus of $A$.
Conversely, let $w = \sum_{i=0}^{n-1} w_iz^i$ be an element of the left nucleus and suppose that $w_k \neq 0$ for some $k$ which is not a multiple of $s$. Then we have $\sigma^k(a) \neq a$. The $z^k$ term of the associator $[w,z^{n-k},z^{k}]$ will be 
\begin{align*}
((w_kz^k)z^{n-k})z^{k} - w_k z^k(z^k z^{n-k}) &=  w_k az^{k} - w_kz^k(a)\\
&= w_kaz^k - w_k \sigma^k(a)z^k
\end{align*}
which is not zero since $\sigma^k(a) \neq a$. Thus the left nucleus of $A$ contains only elements of $B$.
\end{proof}
\end{proposition}

\begin{corollary}
Let $F$ be a finite field and let $L$ be a finite extension of $F$ of degree $n$. The algebra $A = (L/F, \sigma, a)$ is a finite semifield if and only if $a$ belongs to no proper subfield of $L$.
\begin{proof}
If $a$ does not belong to any proper subfield of $L$ then clearly $1, a, a^2, \ldots, a^{n-1}$ are linearly independent over $F$ so Theorem \ref{div} applies. On the other hand, suppose $A$ is a finite semifield and that $a$ belongs to a proper subfield of $L$. Then, using the terminology of the previous proposition, 
\[Nuc_l(A) = L \oplus Lz^s \oplus Lz^{2s} \oplus \cdots \oplus Lz^{n-s}.\] 
This is an associative subalgebra of $A$ so, by Wedderburn's Theorem, this should be a finite field. However, it is easy to check that this subalgebra is not commutative, a contradiction.
\end{proof}
\end{corollary}

We can also determine when two nonassociative cyclic algebras are isomorphic. 

\begin{theorem}\label{iso}
\item Two nonassociative cyclic algebras, $A= (L/F, \sigma, a)$ and $B=(L/F, \sigma, b)$ are isomorphic as $F$-algebras if and only if $\sigma^i(a) = N(l)b$ for some $0 \leq i < n$ and some $l \in L^{\times}$. These isomorphisms are given by
    \[\sum_{j=0}^{n-1} x_j z^j \mapsto \sum_{j=0}^{n-1} \sigma^i(x_j) l\ldots \sigma^{j-1}(l) u^j.\]

\end{theorem}

\begin{corollary}\label{aut}
Let $A = (L/F, \sigma, a)$ be a nonassociative cyclic algebra of degree $n$ and let $l \in L$ be such that $N_{L/F}(l) = 1$. Then every map of the form
\[\sum_{i=0}^{n-1} x_i z^i \mapsto  \sum_{i=0}^{n-1} x_i l \sigma(l) \ldots \sigma^{i-1}(l) z^i\]
is an automorphism of $A$. These maps are the only automorphisms of $A$ unless there exists an element $l' \in L$ such that $\sigma^j(a) = N(l')a$ for some $j = 1, \ldots, n-1$ in which case the map
\[\sum_{i=0}^{n-1} x_i z^i \mapsto  \sum_{i=0}^{n-1} \sigma^j(x_i) l' \sigma(l') \ldots \sigma^{i-1}(l') z^i\]
is also an automorphism of $A$.
\end{corollary}

\section{Nonassociative Cyclic Algebras over Finite Fields}

In this Section we consider the construction of nonassociative cyclic algebras exclusively over finite fields. Given a finite field $F = \mathbb{F}_q$ and a finite extension of degree $n$, $L = \mathbb{F}_{q^n}$, we know we can construct a finite semifield $(L/F, \sigma, a)$ simply by choosing $a$ so that it does not belong to a proper subfield of $L$. For a given field extension $L/F$, we consider two questions: how many nonisomorphic semifields can we construct as a nonassociative cyclic algebra $(L/F, \sigma, a)$ and what is the automorphism group of such a semifield? We give a complete answer of these for the case where $L/F$ is an extension of prime degree.

\begin{remark}
There may be many non-isomorphic finite semifields of order $q^{n^2}$. For example if $n$ is a prime number then we can choose any element $a$ in the field extension and we will get a division algebra. Also  two semifields of the same order can have very different structures, for example consider $q = 2$ and $r = 4$ in the above theorem. Then we can form a semifield with $q^{r^2}= 2^{16}$ elements which has nucleus $\mathbb{F}_{2^4}$ and center $\mathbb{F}_2$. We can also construct a finite semifield with $2^{16}$ elements by letting $q=16$ and $r=2$ in the above theorem. In this case the nucleus will be $\mathbb{F}_{2^8}$ and the center will be $\mathbb{F}_{16}$.
\end{remark}

\begin{example}
Let $F= \mathbb{F}_2$ and let $L = \mathbb{F}_4$, then we can write
\[L = \{0,1,T, 1+T \}\]
where $T^2+T+1 = 0$. Then for the algebra $A_a$ we can either choose $a = T$ or $a = 1+T$, both choices will give a division algebra since $L$ is a field extension of prime degree. We also know that $A_a \cong A_b$ if and only if $\sigma(a) = N(l)b$, but since we are working with finite fields the norm map
\[N: L^{\times} \rightarrow F^{\times}\]
is surjective, so $N(l) = 1$ for all $l \in L^\times$. The statement then reduces to $A_a \cong A_b$ if and only if $\sigma(a) = b$. Now
\[\sigma(T) = T^2 = 1+T.\]
Therefore $A_T \cong A_{1+T}$ so there is only one Sandler semifield up to isomorphism which can be constructed using $L/F$.
\end{example}

The fact that the norm map is surjective for finite field extensions of finite fields allows us to restate the condition from Proposition \ref{iso}: $A_a \cong A_b$ iff $\sigma^i(a) = kb$ for some $0\leq i \leq r-1$ and some $k \in F^{\times}$. We recall the following well-known Lemma.

\begin{lemma}\label{unity}
Let $F = \mathbb{F}_q$ and let $L$ be an extension of $F$ of degree $r$ where $q$ is a prime power and $r$ is prime. $F$ contains a primitive $r$th root of unity if and only if $r$ divides $q-1$.
\end{lemma}

It is well known that if $F$ contains a primitive $r$th root of unity and $L$ is a cyclic field extension of $F$ of degree $r$ (where $r$ is prime to the characteristic of $F$) then $L = F(\omega)$ where $\omega$ is a root of the irreducible polynomial $x^r - a$ for some $a \in F^{\times}$.

\begin{lemma}\label{egnval}
Let $r$ be a prime number and let $F$ be a field of characteristic not $r$ such that $F$ contains a primitive $r$th root of unity. Let $L = F(\omega)$ be a cyclic field extension of $F$ with Galois group $\langle \sigma \rangle$. The eigenvalues of the automorphisms $\sigma^i$ are precisely the $r$th roots of unity. Moreover, the only possible eigenvectors are scalar multiples of the elements $\omega^i$ for $0 \leq i \leq r-1$.
\begin{proof}
Let the elements $1, \omega, \ldots, \omega^{r-1}$ be a basis for $L/F$. The action of $\sigma$ on $\omega$ is given by 
\[\sigma(\omega) = \zeta \omega,\]
where $\zeta$ is a primitive $r$th root of unity. Hence
\[\sigma^i(\omega) = \zeta^i \omega,\]
and 
\[\sigma^i(\omega^j)= \sigma^i(\omega)^j = \zeta^{ij}\omega^j,\]
for all $0 \leq i,j, \leq r-1$. So the $r$th roots of unity are indeed eigenvalues for the automorphisms $\sigma^i$ with eigenvectors $\omega^j$. Now suppose that $k$ is another eigenvalue for $\sigma^i$, i.e. $\sigma^i(x) = kx$ for some $x \in L^{\times}$. Applying the norm map to both sides gives
\[N_{L/F}(\sigma^i(x)) = N_{L/F}(x) = N_{L/F}(kx) = k^rN_{L/F}(x),\]
but this implies that $k^r = 1$ so $k$ is an $r$th root of unity and $k = \zeta^j$ for some $0\leq j \leq r-1$. With our chosen basis of $L/F$ the matrix of the automorphism $\sigma^i$ is
\[\left( \begin{array}{ccccc}
1 & 0 & 0 & \cdots & 0 \\
0 & \zeta^i & 0 & \cdots & 0 \\
0 & 0 & \zeta^{2i} & \cdots & 0\\
\vdots & \vdots & \vdots & \ddots & \vdots \\
0 & 0 & 0 & \cdots & \zeta^{(r-1)i} \end{array} \right). \]
The equation $\sigma^i(x) = \zeta^j x$ in matrix form becomes:
\[\left( \begin{array}{ccccc}
1 & 0 & 0 & \cdots & 0 \\
0 & \zeta^i & 0 & \cdots & 0 \\
0 & 0 & \zeta^{2i} & \cdots & 0\\
\vdots & \vdots & \vdots & \ddots & \vdots \\
0 & 0 & 0 & \cdots & \zeta^{(r-1)i} \end{array} \right) \begin{pmatrix} x_0 \\ x_1 \\ x_2 \\ \vdots \\ x_{r-1} \end{pmatrix} = \zeta^j \Big(x_0, x_1, x_2, \ldots, x_{r-1}\Big), \]
where $x_i \in F$ and $x = (x_0, x_1, \ldots, x_{r-1})$ is written as an $r$-tuple with respect to our basis. This gives
\[\Big(x_0, \zeta^i x_1, \zeta^{2i}x_2, \ldots, \zeta^{(r-1)i} x_{r-1} \Big) = \Big(\zeta^j x_0 , \zeta^j x_1, \zeta^j x_2, \ldots, \zeta^j x_{r-1} \Big),\]
which implies that all the $x_k$ are zero except for one, say $x_{k_0}$, where $k_0$ is such that $k_0 i \cong j$ mod $r$. Hence $x = x_{k_0} \omega^{k_0}$ as required.
\end{proof}
\end{lemma}

In the case of the semifields defined above we are looking for elements $a \in L \setminus F$ with $\sigma^i(a) = k a$. Note that $k=1$ is not relevant for this case since $\sigma^i(a) = a$ if and only if $a \in F$. Define an equivalence relation on the set $L \setminus F$ by
\[a \sim b \text{ if and only if } (L/F, \sigma, a) \cong (L/F, \sigma, b).\]
We wish to know how many distinct equivalence classes there are in $L \setminus F$ for a given finite field $F$ and extension $L$ of prime degree.
\begin{theorem}
Let $F = \mathbb{F}_q$ and let $L$ be an extension of $F$ of degree $r$ where $r$ is prime and $q$ is a prime power. If $r$ divides $q-1$ then there are 
\[r-1 + \frac{q^r-q - (q-1)(r-1)}{r(q-1)}\]
equivalence classes of the above equivalence relation. If $r$ does not divide $q-1$ then there are precisely 
\[\frac{q^r-q}{r(q-1)}\]
equivalence classes.
\begin{proof}
For each $a \in L \setminus F$ we have 
\[(L/F, \sigma, a) \cong (L/F, \sigma, \sigma^i(a))\]
for $0 \leq i \leq r-1$ and 
\[(L/F, \sigma, a) \cong (L/F, \sigma, ka)\]
for $k \in F^{\times}$. If the elements $k\sigma^i(a)$ for $0 \leq i \leq r-1$ and $k \in F^{\times}$ are all distinct then there are precisely $r(q-1)$ elements in the equivalence class of $a$. If they are not all distinct then this is equivalent to $\sigma^i(a) = ka$ for some $i$ and some $k \in F^{\times}$. We saw in the proof of Lemma \ref{egnval} that if $\sigma^i(a) = ka$ then $k$ is an $r$th root of unity. As mentioned above we cannot have $k = 1$ so $k$ is a primitive $r$th root of unity. From Lemma \ref{unity} we know this happens if and only if $r$ divides $q-1$. Hence if $r$ does not divide $q-1$ then from $q^r-q$ elements in $L \setminus F$ we get 
\[\frac{q^r - q}{r(q-1)}\]
equivalence classes. On the other hand, if $r$ does divide $q-1$ then $F$ contains the primitive $r$th roots of unity and so we may write $L$ as 
\[F[T]/(T^r - c)\]
for some $c \in F^{\times}$. Lemma \ref{unity} tells us that the only elements $a \in L \setminus F$ with $\sigma^i(a) = ka$ are the elements $T^j$ for $1 \leq j \leq r-1$ and scalar multiples of these. Moreover for each $T^j$ we have $\sigma^i(T^j) = \zeta^{ij}T^j$ and $\zeta^{ij} \in F$ so there are only $q-1$ distinct elements in the equivalence class of each $T^j$. Hence the $(q-1)(r-1)$ elements $kT^j$ ($k \in F^{\times}$ and $1 \leq j \leq r-1$) form exactly $r-1$ equivalence classes. Since these are all the elements in $L \setminus F$ which are eigenvectors for the automorphisms $\sigma^i$ we can deduce that remaining $q^r-q - (q-1)(r-1)$ elements will form 
\[\frac{q^r-q - (q-1)(r-1)}{r(q-1)}\]
equivalence classes. In total we obtain
\[r-1 + \frac{q^r-q - (q-1)(r-1)}{r(q-1)}\]
equivalence classes.
\end{proof}
\end{theorem}

\begin{corollary} \label{numb}
Let $F = \mathbb{F}_q$ and let $L$ be an extension of $F$ of degree $r$ where $r$ is prime and $q$ is a prime power. If $r$ divides $q-1$ then there are 
\[r-1 + \frac{q^r-q - (q-1)(r-1)}{r(q-1)}\]
non-isomorphic semifields arising from the construction $(L/F, \sigma, a)$ . If $r$ does not divide $q-1$ then there are precisely 
\[\frac{q^r-q}{r(q-1)}\]
non-isomorphic semifields from this construction.
\end{corollary}

We now move on to the question of determining the automorphism group of a given semifield $(L/F, \sigma, a)$. Recall from Corollary \ref{aut} that automorphisms of $(L/F, \sigma, a)$ are given by
\[\sum_{i=0}^{n-1} x_i z^i \mapsto  \sum_{i=0}^{n-1} x_i l \sigma(l) \ldots \sigma^{i-1}(l) z^i\]
where for some $l \in L$ with $N_{L/F}(l) = 1$. These are all the automorphisms unless there exists an $l' \in L$ such that $\sigma^i(a) = N_{L/F}(l') a$ in which case 
\[\sum_{i=0}^{n-1} x_i z^i \mapsto  \sum_{i=0}^{n-1} \sigma^i(x_i) l' \sigma(l') \ldots \sigma^{i-1}(l') z^i\]
is also an automorphism. Hence the automorphisms of $(L/F, \sigma, a)$ also depend on the existence of elements $k \in F^{\times}$ such that $\sigma^i(a) = ka$. However, it is clear that the kernel of the norm map $N_{L/F}$ will be isomorphic to a subgroup of the automorphism group of the semifield. We recall the following well-known fact.

\begin{proposition}
Let $F = \mathbb{F}_q$ and let $L$ be an extension of $F$ of degree $n$: $L = \mathbb{F}_{p^n}$. The kernel of the norm map $N_{L/F}$ is a cyclic subgroup of order $s = \frac{q^n-1}{q-1}$.
\end{proposition}

It follows from the multiplicativity of the norm that for every $k \in F^{\times}$ there exist exactly $\frac{q^n-1}{q-1}$ elements $x \in L$ with $N_{L/F}(x) = k$.

\begin{corollary} 
Let $F = \mathbb{F}_q$ and let $L$ be an extension of $F$ of prime degree $r$. If $r$ does not divide $q-1$ then for all $a \in L \setminus F$, Aut$((L/F, \sigma, a)) \cong \mathbb{Z}/s \mathbb{Z}$ where $s = \frac{q^r-1}{q-1}$.
\begin{proof}
All automorphisms of $(L/F, \sigma, a)$ correspond to elements of ker$(N_{L/F})$ unless there exists $k \in F^{\times}$ with $\sigma^i(a) = k a$ for some $1 \leq i \leq r-1$. It was shown in the proof of Lemma \ref{egnval} that any such $k$ is a primitive $r$th root of unity which cannot happen in $F$ by Lemma \ref{unity}. The result now follows from the previous proposition.
\end{proof}
\end{corollary}

Consider now the case where $F$ and $L$ are as above but $r$ does divide $q-1$. Since $F$ contains all $r$th roots of unity we may write 
\[L = F[T]/(T^r - c)\]
for some $c \in F^{\times}$. Define the set
\[S:= \{T^j \mid 1 \leq j \leq r-1\} \subset L.\]

\begin{corollary}\label{auto1}
Let $F$ and $L$ be as above with
\[L = F[T]/(T^r - c).\]
Then for all $a \in L \setminus (F \cup S)$, Aut$((L/F, \sigma, a)) \cong \mathbb{Z}/s \mathbb{Z}$ where $s = \frac{q^r-1}{q-1}$. If $a \in S$ then Aut$((L/F, \sigma, a))$ is a group of order 
\[r \frac{q^r-1}{q-1}.\]
\begin{proof}
Lemma \ref{egnval} states that the only elements $a \in L \setminus F$ with $\sigma^i(a) = ka$ for some $1 \leq i \leq r-1$ and some $k \in F$ are the elements of $S$, hence the first claim follows. Now suppose $a = T^j$ for some $j$. For each $i$ from $0$ to $r-1$ there exists a unique  $r$th root of unity $\zeta_i$ such that $\sigma^i(a) = \zeta_i a$ (note $\zeta_0 = 1$). There are exactly $\frac{q^r-1}{q-1}$ elements $l \in L$ with $N_{L/F}(l) = \zeta_i$ and each of these elements correspond to a unique automorphism. Therefore in total we have
\[r \frac{q^r-1}{q-1}\]
automorphisms.
\end{proof}
\end{corollary}

\begin{example}
Let $F = \mathbb{F}_3$ and $L = \mathbb{F}_9$ where
\[L = F[T]/(T^2 - 2) = \{0,1,2,T,2T, T+1, T+2, 2T+1, 2T+2\}.\]
There are two nonisomorphic Sandler semifields for $L/F$. These are\\ $A_T := (L/F, \sigma, T)$ and $A_{T+1}:= (L/F, \sigma, T+1)$. Moreover, Aut($A_{T+1}$) is the cyclic group of order $4$ whereas Aut($A_T$) is isomorphic to the group  of quaternion units.
\end{example}
\begin{proof}
Corollary \ref{numb} tells us that there are $2$ nonisomorphic semifields of order $81$ arising from the construction $(L/F, \sigma, a)$ and Lemma \ref{egnval} implies that two such nonisomorphic semifields are $(L/F, \sigma, T)$ and $(L/F, \sigma, T+1)$. Denote these two semifields by $A_T$ and $A_{T+1}$ respectively. Now Corollary \ref{auto1} tells us that the automorphism group of $A_{T+1}$ will be the cyclic group of order $4$. However, the automorphism group of $A_T$ will be a group of order $8$. To calculate what this group is we introduce the following notation:

Let $l \in L$ be such that $N_{L/F}(l)  = 1$. Denote by $\varphi_l$ the automorphism of $A_T$ given by
\[\varphi_l(x_0 + x_1 z) = x_0 + x_1 l z\]
for all $x = x_0 + x_1 z \in A_T$. Now let $m \in L$ be such that $\sigma(T) = N_{L/F}(m)T$. Denote by $\theta_m$ the automorphism of $A_T$ given by
\[\theta_m(x_0 + x_1 z) = \sigma(x_0) + \sigma(x_1) m z.\]
Since $\sigma(T) = T^3 = 2T$ we require all those $m \in L$ with $N_{L/F}(m) = 2$. A quick calculation shows that these are:
\[m \in  \{1+T, 1+2T, 2+T, 1+2T\}.\]
Moreover the set of elements of $L$ with norm $1$ is:
\[\{ 1, 2, T, 2T\}.\]

Hence using the above notation the automorphism group will be
\[Aut(A_T) = \{\varphi_1, \varphi_2, \varphi_T, \varphi_{2T}, \theta_{1+T}, \theta_{1+2T}, \theta_{2+T}, \theta_{2+2T}\}.\]
It is easy to check that this is a non-abelian group of order eight, also it has one element of order two, namely $\varphi_2$ and the rest of the (non-identity) elements are of order four. The classification of small groups then implies that Aut($A_T) \cong \mathcal{Q}$, the group of quaternion units.
\end{proof}

\begin{theorem}
Let $F = \mathbb{F}_q$ be a field of characteristic not $2$ and let $L$ be a quadratic extension of $F$. For $a \in L \setminus F$ put $A_a := (L/F, \sigma, a)$. Then Aut($A_a)$ is the cyclic group of order $(q+1)$ or the dicyclic group of order $2q+2$
\begin{proof}
Write $L = F[T]/(T^2 - c)$ for some $c \in F^{\times}$. If $a = kT$ for some nonzero $k \in F$ then we know from Corollary \ref{auto1} that Aut($A_a)$ is of order $2(q+1)$, otherwise Aut($A_a$) is the cyclic group of order $q+1$. We may assume that $a = T$ since \newline
$A_a \cong A_{ka}$ for all nonzero $k \in F$. Since $\sigma(T) = -T$ we have the following automorphisms:
\[\varphi_{l_i}:A_a \rightarrow A_a: \hspace{15mm} x_0 + x_1z \mapsto x_0 + x_1 l_i z,\]
for $x_0 + x_1 z \in A_a$ and $l_i \in L$ such that $N_{L/F}(l_i) = 1$. There are precisely $q+1$ such maps. We also have the automorphisms 
\[\theta_{m_j}:A_a \rightarrow A_a: \hspace*{15mm} x_0 + x_1 z \mapsto \sigma(x_0) + \sigma(x_1) m_i z\]
for all $m_j \in L$ such that $N_{L/F}(m_j) = -1$. Again there are precisely $q+1$ such maps. We note the following relations between the automorphisms:
\[\varphi_{l_i} \circ \varphi_{l_j} = \varphi_{l_il_j}, \hspace*{20mm} \varphi_{l_i} \circ \theta_{m_j} = \theta_{l_i m_j},\]
\[\theta_{m_i} \circ \varphi_{l_j}, = \theta_{m_i \sigma (l_j)} \hspace*{20mm} \theta_{m_i} \circ \theta_{m_j} = \varphi_{m_i \sigma(m_j)}.\]
Recall that we can describe the dicyclic group of order $4n$, $Dic_n$, by the following presentation:
\[ Dic_n = \langle x, y \mid y^{2n} = 1, x^2 = y^n, x^{-1}yx = y^{-1} \rangle.\]
We claim that Aut($A_a) \cong Dic_n$ for $n = (q+1)/2$. First note that the group Ker$(N_{L/F}) = \{l_i \mid N_{L/F}(l_i) = 1\}$ is cyclic, so pick $l_0 \in $ Ker$(N_{L/F})$ which generates it as a group. Also pick any $m_0$ such that $N_{L/F}(m_0) = -1$, then the map 
\[\varphi_{l_0} \mapsto y, \hspace*{20mm} \theta_{m_0} \mapsto x\]
is an isomorphism from $Aut(A_a) \rightarrow Dic_q$. Clearly we have
\[ (\varphi_{l_0})^{2n} = \varphi_{(l_0)^{2n}} = \varphi_1 = Id_{A_a}.\]
From this it follows that $(\varphi_{l_0})^n = \varphi_{-1}$ and so
\[(\theta_{m_0})^2 = \varphi_{N_{L/F}(m_0)} = \varphi_{-1} = (\varphi_{l_0})^n.\]
Finally, note that $(\theta_{m_0})^{-1} = \theta_{m_j}$ where $m_j$ is such that $m_j \sigma(m_0) = 1$ and $(\varphi_{l_0})^{-1} = \varphi_{\sigma(l_0)}$. Hence
\[(\theta_{m_0})^{-1} \circ \varphi_{l_0} \circ \theta_{m_0} = \varphi_{m_j \sigma(l_0) \sigma(m_0)} = \varphi_{\sigma(l_0)} = (\varphi_{l_0})^{-1}.\]
\end{proof}
\end{theorem}

\section{Hughes-Kleinfeld and Knuth Semifields}
In \cite{MR0117261}, Hughes and Kleinfeld give a construction of a finite semifield which is quadratic over a finite field $L$ contained in the right and middle nucleus. In his thesis \cite{MR2939390} (see also \cite{MR0175942}), Knuth considered three similar constructions of finite semifields. These three along with the Hughes-Kleinfeld semifields are some of the best known constructions of semifields and have been studied in a variety of contexts (\cite{MR2336394}, \cite{MR2942770}, for example). Under certain conditions they each possess different combinations of left, right and middle nuclei and so they are mutually non-isomorphic. In this Section we give the four constructions over arbitrary (possibly infinite) fields. 

Let $F$ be a field and let $L$ be a separable field extension of $F$. We consider four multiplications on the $F$-vector space 
\[L \oplus L.\]
Pick elements $\eta$ and $\mu$ of $L$ and a nontrivial automorphism $\sigma \in Aut(L/F)$. For elements $x,y,u,v \in L$ the four multiplications are given as follows:
\begin{align*}
Kn_1: (x,y) \circ (u,v) &= (xu + \eta \sigma(v)\sigma^{-2}(y) , vx + y \sigma(u) + \mu \sigma(v)\sigma^{-1}(y) ),\\ \\
Kn_2: (x,y) \circ (u,v) &= (xu + \eta \sigma^{-1}(v) \sigma^{-2}(y) , vx + y \sigma(u) + \mu v\sigma^{-1}(y)),\\ \\
Kn_3: (x,y) \circ (u,v) &= (xu + \eta \sigma^{-1}(v)y , vx + y \sigma(u) + \mu vy),\\ \\
HK : (x,y) \circ (u,v) &= (xu + \eta \sigma(v)y , vx + y \sigma(u) + \mu \sigma(v)y ).
\end{align*}
We will denote the vector-space $L \oplus L$ endowed with each of the above multiplications by $Kn_1(L, \sigma, \eta, \mu), Kn_2(L, \sigma, \eta, \mu), Kn_3(L, \sigma, \eta, \mu)$ and $HK(L, \sigma, \eta, \mu)$, respectively. This notation reflects the fact that the first three are the constructions defined by Knuth and the last one is the construction defined by Hughes-Kleinfeld. If it is clear from the context or irrelevant to the discussion we may omit some or all of the parameters. Each of $Kn_1, Kn_2, Kn_3, HK$ are unital $F$-algebras with unit element $(1,0)$. Each of them also contain $L \oplus 0$ as a subalgebra which we will identify with the field $L$.

\begin{theorem}[\cite{MR0117261}, \cite{MR0175942}]
Each of the above constructions is a division algebra if and only if the equation
\[w \sigma(w) + \mu w - \eta\]
has no solutions in $L$.
\end{theorem}
Thus we get four constructions for division algebras which will be finite semifields if $F$ is finite. 

By the previous theorem, if $\eta = 0$ then $Kn_1, Kn_2, Kn_3$ and $HK$ are never semifields for any choice of $\mu$. Hence we shall assume that $\eta \neq 0$ from now on.

\section{Nuclei}
In this section we calculate some of the nuclei of the algebras $Kn_1, Kn_2, Kn_3$ and $HK$ and show that no two of them possess the same combination of left, right and middle nuclei unless $\sigma^2 = Id$ and $\mu = 0$. If both $\sigma^2 = Id$ and $\mu = 0$ then the multiplication for each algebra is the same:
\[(x,y) \circ (u,v) = (xu + \eta y \sigma(v), xv + y \sigma(u)).\]
In this case $\sigma$ is of order two in $Gal(L/F)$ and the Fundamental Theorem of Galois Theory tells us that $L$ has a subfield $E$ such that $[L:E] = 2$ and $Gal(L/E) = \{Id, \sigma\}$. Hence, the multiplication given above defines a quaternion algebra over $E$ which is associative if $\eta \in E$ and nonassociative if $\eta \in L \setminus E$. The nonassociative case is covered in the previous sections (they are nonassociative cyclic algebras of degree $2$). The theory of associative quaternion algebras is well known.

\begin{proposition} Suppose that either $\sigma^2 \neq Id$ or that $\mu \neq 0$.\\
(i) $L$ is equal to the middle and right nucleus of $Kn_2(L, \sigma, \eta, \mu)$ but is not contained in the left nucleus.\\
(ii) $L$ is equal to the left and right nucleus of $Kn_3(L, \sigma, \eta, \mu)$ but is not contained in the middle nucleus.\\
(iii) $L$ is equal to the left and middle nucleus of $HK(L, \sigma, \eta, \mu)$ but is not contained in the right nucleus.\\
(iv) $L$ is not contained in the left, right or middle nucleus of $Kn_1(L, \sigma, \eta, \mu)$.
\begin{proof}
This result has been proved in various forms for finite fields in the above mentioned papers. Here we give a proof of (i) for completeness since it is crucial to the results of the next section. The proofs of (ii), (iii) and (iv) are similar. We will also drop the circle notation for multiplication and denote it simply by juxtaposition: $(x,y) \circ (u,v)= (x,y)(u,v)$. To show $L$ is contained in the middle nucleus we calculate
\begin{align*}
((x,y)(l,0))(u,v) &= (xl,y \sigma(l))(u,v)\\
&= (xlu+ \eta \sigma^{-2}(y)\sigma^{-1}(l) \sigma^{-1}(v), xlv + y \sigma(l)\sigma(u) + \mu \sigma^{-1}(y)lv).
\end{align*}
On the other hand
\begin{align*}
(x,y)((l,0)(u,v)) &= (x,y)(lu,lv)\\
&= (xlu+ \eta \sigma^{-2}(y)\sigma^{-1}(lv) , xlv + y \sigma(lu) + \mu \sigma^{-1}(y)lv).
\end{align*}
These two expressions are the same hence $L \subseteq Nuc_m(Kn_2)$. To show that there are no other elements in the middle nucleus of $Kn_2$ it suffices to check that no elements of the form $(0,m)$ belong to the middle nucleus. This is because the associator is linear:
\[[(x,y), (l,m), (u,v)] = [(x,y), (l,0), (u,v)] + [(x,y), (0,m), (u,v)].\]
If $(0,m) \in Nuc_m(Kn_2)$ then for all $v \in L$ we should have
\[[(0,v), (0,m), (0,1)] = (0,0),\]
however, calculating directly we get
\begin{align*}
&((0,v)(0,m))(0,1)  - (0,v)((0,m)(0,1)) =\\
&(\eta \sigma^{-2}(\mu) \sigma^{-3}(v)\sigma^{-2}(m), \eta \sigma^{-2}(v) \sigma^{-2}(m) + \mu \sigma^{-1}(\mu) \sigma^{-2}(v)\sigma^{-1}(m))\\
&- (\eta \sigma^{-2}(v) \sigma^{-1}(\mu) \sigma^{-2}(m), v \sigma(\eta)\sigma^{-1}(m) + \mu \sigma^{-1}(v) \mu \sigma^{-1}(m)).
\end{align*}
If $\mu \neq 0$ then looking at the first term gives
\[\sigma^{2}(v) \sigma^{-1}(\mu) = \sigma^{-2}(\mu) \sigma^{-3}(v)\]
for all $v \in L$ which can't hold. If $\mu = 0$ then by hypothesis we must have $\sigma^2 \neq Id$ and looking at the second term gives the equation
\[\eta \sigma^{2}(v) = v \sigma(\eta)\]
which again can't hold for all $v \in L$. A similar calculation shows that the right nucleus is also equal to $L$. To show that $L$ is not contained in the left nucleus we check
\begin{align*}
((l,0)(x,y)))(u,v) &= (lx,ly)(u,v)\\
&= (lxu+ \eta \sigma^{-2}(ly) \sigma^{-1}(v), lxv + ly \sigma(u) + \mu \sigma^{-1}(ly)v),
\end{align*}
whereas
\begin{align*}
(l,0)((x,y)(u,v)) &= (l,0)(xu + \eta \sigma^{-2}(y) \sigma^{-1}(v), xv + y \sigma(u) + \mu \sigma^{-1}(y) v)\\
&= (lxu+ l\eta \sigma^{-2}(y) \sigma^{-1}(v), lxv + ly \sigma(u) + l\mu \sigma^{-1}(y)v).
\end{align*}
These equations are not equal for all $l \in L$ unless $\sigma^2 = Id$ and $\mu = 0$.
\end{proof}
\end{proposition}

\begin{remark} \label{opp}
In their paper, \cite{MR0117261}, Hughes and Kleinfeld show that the right and middle nuclei of their algebra are equal to $L$. This is because they use a slightly different definition of multiplication to us. They consider the product
\[(x,y) \circ (u,v) = (xu + \eta \sigma(y) v, yu + \sigma(x) v + \mu \sigma(y) v),\]
for all $x,y,u,v \in L$. It is easily checked that this multiplication gives the opposite algebra $HK^{op}$ which explains the swap of right and left nucleus.
\end{remark}

\begin{corollary}
The algebras $Kn_1, Kn_2, Kn_3, HK$ are mutually non-isomorphic unless $\sigma^2 = Id$ and $\mu = 0$, in which case they are the same algebra.
\begin{proof}
Since isomorphisms must preserve each of the left, right and middle nuclei, the first claim holds when either $\sigma^2 \neq Id$ or $\mu \neq 0$. On the other hand if both $\sigma^2 = Id$ and $\mu = 0$ the the definition of the multiplication in each semifield is exactly the same.
\end{proof}
\end{corollary}

\section{Automorphisms}
In this section we describe all automorphisms for the algebras $HK, Kn_2$ and $Kn_3$. We also exhibit some automorphisms for the algebra $Kn_1$.

\begin{proposition}\label{hkauto}
Let $A = HK(L, \sigma, \eta, \mu)$. All automorphisms of $A$ are of the form 
\[(x,y) \mapsto (\tau(x), \tau(y)b),\]
where $\tau \in Aut(L/F)$ commutes with $\sigma$ and $b \in L^\times$ is such that
\[ \eta b \sigma(b) = \tau(\eta) \text{ and } \mu \sigma(b) = \tau(\mu).\]
\begin{proof}
Let $\varphi:A \rightarrow A$ be an automorphism. Since $Nuc_l(A) = L$ and isomorphisms preserve left, right and middle nuclei, we must have $\varphi(L) = L$. Hence $\varphi|_L = \tau \in Aut(L/F)$. We can write any element $(x,y)$ of $HK$ as 
\[(x,y) = (x,0) + (y,0)(0,1).\]
Since $\varphi((x,0)) = (\tau(x), 0)$ for all $x \in L$, it remains to determine $\varphi((0,1))$. Suppose $\varphi((0,1)) = (a,b)$ for some $a,b \in L$, thus we can write
\begin{align*}
\varphi((x,y)) &= \varphi((x,0)) + \varphi((y, 0))\varphi((0,1))\\
&= (\tau(x), 0) + (\tau(y), 0)(a,b)\\
&= (\tau(x) + \tau(y)a, \tau(y)b).
\end{align*}
For all $m \in L$ we must have
\[\varphi((0,1)(m,0)) = \varphi((0,1))\varphi((m,0)).\]
On the one hand we have 
\[\varphi((0,1)(m,0)) = \varphi((0,\sigma(m)) = (\tau(\sigma(m))a, \tau(\sigma(m))b).\]
Whereas the right hand side becomes 
\[(a,b)(\tau(m), 0) = (\tau(m)a, \sigma(\tau(m))b).\]
Since $\sigma \neq Id$, this implies that $a = 0$ and that $\tau$ and $\sigma$ commute. Finally, we have $(0,1)(0,1) = (\eta, \mu)$ and so $\varphi((0,1)(0,1)) = (\tau(\eta), \tau(\mu)b)$. However,
\[\varphi((0,1))\varphi((0,1)) = (0,b)(0,b) = (\eta b\sigma(b), \mu b \sigma(b)),\]
and so we arrive at the conditions $\eta b \sigma(b) = \tau(\eta)$ and $\mu \sigma(b) = \tau(\mu)$.

Conversely, it is not difficult to check that the maps given above are indeed automorphisms.

\end{proof}
\end{proposition}

\begin{proposition}\label{kn2auto}
Let $A = Kn_2(L, \sigma, \eta, \mu)$. All automorphisms of $A$ are of the form 
\[(x,y) \mapsto (\tau(x), \tau(y)b),\]
where $\tau \in Aut(L/F)$ commutes with $\sigma$ and $b \in L^{\times}$ is such that
\[ \eta b \sigma^{-1}(b) = \tau(\eta) \text{ and } \mu b = \tau(\mu).\]
\end{proposition}

\begin{proposition}\label{kn3auto}
Let $A = Kn_3(L, \sigma, \eta, \mu)$. All automorphisms of $A$ are of the form 
\[(x,y) \mapsto (\tau(x), \tau(y)b),\]
where $\tau \in Aut(L/F)$ commutes with $\sigma$ and $b \in L^{\times}$ is such that
\[ \eta \sigma^{-1}(b) \sigma^{-2}(b) = \tau(\eta) \text{ and } \mu \sigma^{-1}(b) = \tau(\mu).\]
\end{proposition}

The proof of these propositions follow the same line of argument as Proposition \ref{hkauto} with the only difference coming at the end when we deduce what conditions the element $b \in L$ must satisfy. The key part in these proofs is using the fact that either the left, right or middle nucleus of $HK, Kn_2$ and $Kn_3$ is equal to $L$. From this we deduce that any automorphism of the algebra must restrict to an automorphism on $L$. For $Kn_1$, $L$ is not contained in any of the nuclei so we cannot make this deduction. However, if we assume that an automorphism of $Kn_1$ restricts to an automorphism of $L$, then it must be of a similar form to the above maps.

\begin{proposition}
Let $A = Kn_1(L, \sigma, \eta, \mu)$ and suppose $\varphi$ is an automorphism of $A$ which restricts to an automorphism of $L$: $\varphi|_L =  \tau \in Aut(L/F)$. Then $\tau$ commutes with $\sigma$ and for all $(x,y) \in A$ 
\[\varphi((x,y)) = (\tau(x), \tau(y)b),\]
where $\eta \sigma^{-1}(b) \sigma^{-2}(b) = \tau(\eta)$ and $\mu \sigma(b) \sigma^{-1}(b) = \tau(\mu) b$.
\end{proposition}

The proof of this is also similar to that of Proposition \ref{hkauto}. It is not yet clear if these are all automorphisms of $Kn_1$.

Whenever $\mu \neq 0$, there will be very few automorphisms for each of these algebras, in many cases the automorphism group of the algebra will be smaller than the Galois group of $L/F$. The exact size of the automorphism group depends on the position of the elements $\eta$ and $\mu$ within $L$.

\begin{proposition}
Let $G = Aut(L/F)$ and let $C_G(\sigma)$ be the centraliser of $\sigma$ in $G$. If $A$ is one of the algebras $HK(L, \sigma, \eta, \mu), Kn_2(L, \sigma, \eta, \mu)$ or $Kn_3(L, \sigma, \eta, \mu)$ where $\mu \neq 0$, then the automorphism group of $A$ is isomorphic to the subgroup of $C_G(\sigma)$ which fixes the element $\mu \sigma(\mu) \sigma(\eta)^{-1}$, i.e.
\[Aut(A) \cong \{\tau \in C_G(\sigma) \mid \tau \Bigg( \frac{\mu \sigma(\mu)}{\sigma(\eta)}\Bigg) = \frac{\mu \sigma(\mu)}{\sigma(\eta)}\}.\]
\begin{proof}
Suppose $A = HK(L, \sigma, \eta, \mu)$ and denote by $\varphi_\tau^b$ the automorphism of $A$ 
\[(x,y) \mapsto (\tau(x), \tau(y)b),\]
for all $(x,y) \in A$. From Proposition \ref{hkauto} we know that $\tau \in C_G(\sigma)$ and $\mu \sigma(b) = \tau(\mu)$. Since $\mu \neq 0$, the element $b \in L$ is determined completely by the action of $\tau$ on $\mu$ so we may drop the superscript $b$ in $\varphi_\tau^b$ and write $\varphi_\tau$. We also have the equation $\eta b \sigma(b) = \tau(\eta)$ defining $\varphi_\tau$. Substituting in $ b = \sigma^{-1}(\tau(\mu)) \sigma^{-1}(\mu)^{-1}$ and rearranging gives
\[\sigma(\eta) \tau(\mu)\sigma(\tau(\mu)) = \sigma(\tau(\eta))\mu \sigma(\mu).\]
Since $\sigma$ and $\tau$ commute, we can rearrange this further to get
\[\tau \Bigg(\frac{\mu \sigma(\mu)}{\sigma(\eta)}\Bigg) = \frac{\mu \sigma(\mu)}{\sigma(\eta)}.\]
Now it is a straightforward calculation to check that if $\varphi_{\tau_1}$ and $\varphi_{\tau_2}$ are two such automorphisms defined above then,
\[\varphi_{\tau_1} \circ \varphi_{\tau_2} = \varphi_{\tau_1 \tau_2}.\]
Hence, $\varphi_\tau \mapsto \tau$ is the required isomorphism.
\end{proof}
\end{proposition}

\begin{corollary}\label{quad}
Let $L/F$ be a quadratic, separable field extension and suppose $A$ is one of the algebras $HK(L, \sigma, \eta, \mu), Kn_2(L, \sigma, \eta, \mu)$ or $Kn_3(L, \sigma, \eta, \mu)$ where $\mu \neq 0$, then the automorphism group of $A$ is isomorphic to $\mathbb{Z}/2\mathbb{Z}$ if $\eta \in F$. Otherwise, $Aut(A) = \{Id\}$.
\begin{proof}
Since $L/F$ is quadratic, $\sigma$ is the nontrivial automorphism of $L/F$. By the previous proposition, we can have at most two possible automorphisms of $A$: $\varphi_{Id}$ and $\varphi_{\sigma}$. Moreover, $\varphi_{\sigma} \in Aut(A)$ if and only if 
\[\sigma \Bigg(\frac{\mu \sigma(\mu)}{\sigma(\eta)}\Bigg) = \frac{\mu \sigma(\mu)}{\sigma(\eta)}.\]
Now $\mu \sigma(\mu) \in F^{\times}$ for all $\mu \in L^{\times}$, so this condition is equivalent to 
\[\sigma \Bigg(\frac{1}{\sigma(\eta)}\Bigg) = \frac{1}{\sigma(\eta)}\]
i.e. $\eta = \sigma(\eta)$. This happens if and only if $\eta \in F$.
\end{proof}
\end{corollary}

\bibliographystyle{alpha}
\bibliography{refs}
\end{document}